\documentclass{gtart_h}  

\def\ifplaintex{\expandafter\ifx\csname documentclass\endcsname\relax}

\ifplaintex 
\else
\fi

\expandafter\ifx\csname beginpicture\endcsname\relax
\expandafter\ifx\csname documentclass\endcsname\relax
\input pictex \else
\input prepictex \input pictex \input postpictex \fi\fi

\def\gt{{\mathsurround=0pt\it $\cal G\mskip-2mu$eometry \&\ 
$\cal T\!\!$opology}}        

\def\gtp{{\mathsurround=0pt\it $\cal G\mskip-2mu$eometry \&\ 
$\cal T\!\!$opology $\cal P\!$ublications}}  

\def\lognumber#1{\def\thelognumber{#1}}
\def\volumenumber#1{\def\thevolumenumber{#1}}
\def\papernumber#1{\def\thepapernumber{#1}}
\def\volumeyear#1{\def\thevolumeyear{#1}}

\def\pagenumbers#1#2{\def\startpage{#1}\def\finishpage{#2}}
\def\published#1{\def\publishdate{#1}}
\def\proposed#1{\def\theproposer{#1}}
\def\seconded#1{\def\theseconders{#1}}
\def\received#1{\def\receiveddate{#1}}
\def\revised#1{\def\reviseddate{#1}}
\def\accepted#1{\def\accepteddate{#1}}

\long\def\asciiabstract#1{\long\def\theasciiabstract{#1}}
\def\asciikeywords#1{\def\theasciikeywords{#1}}

\let\\\par\let\thelognumber\relax
\let\thevolumenumber\relax\let\thepapernumber\relax
\let\thevolumeyear\relax\let\thesamplenumber\relax\let\startpage\relax
\let\finishpage\relax\let\publishdate\relax\let\receiveddate\relax
\let\reviseddate\relax\let\accepteddate\relax\let\theasciititle\relax
\let\theasciiauthors\relax
\let\theasciiabstract\relax\let\theasciikeywords\relax
\let\theasciiemail\relax\let\theshortauthors\relax\let\theshorttitle\relax

\long\def\maketitlep{   

\count0=\startpage

\gt\hfill      
\beginpicture
\setcoordinatesystem units <0.33truein, 0.33truein> point at 2.2 0.9
\setplotsymbol ({$\cal G$})
\plotsymbolspacing=9truept
\circulararc 315 degrees from 0 1 center at 0 0
\setplotsymbol ({$\cal T$})
\circulararc 315 degrees from 1 -1 center at 1 0
\endpicture
\break
{\small\ifx\thesamplenumber\relax 
Volume \else Sample
\fi\thevolumenumber\ (\thevolumeyear)
\startpage--\finishpage\nl
Published: \publishdate}
\vglue 0.5truein plus 0.4fil minus 0.1truein

{\parskip=0pt\leftskip 0pt plus 1fil\def\\{\par\smallskip}{\ifplaintex\large
\else\Large\fi\bf\thetitle}\par\medskip}   

\vglue 0pt plus 0.1fil 

{\parskip=0pt\leftskip 0pt plus 1fil\def\\{\par}{\sc\theauthors}
\par\medskip}

\vglue 0pt plus 0.1fil 

{\small\parskip=0pt\let\newline\\
{\leftskip 0pt plus 1fil\def\\{\par}{\sl\theaddress}\par}
\expandafter\ifx\theemail\relax    
\relax\else\vglue 5pt plus 0.02fil minus 2pt\def\\{\stdspace{\rm 
and}\stdspace} 
\cl{Email:\stdspace\tt\theemail}\fi
\ifx\theurl\relax                  
\relax\else\vglue 5pt plus 0.02fil minus 2pt\def\\{\stdspace{\rm 
and}\stdspace}
\cl{URL:\stdspace\tt\theurl}\fi\par}

\vglue 7pt plus 0.3fil minus 3pt

{\bf Abstract}
\vglue 5pt plus 0.1fil minus 2pt

\theabstract

\vglue 7pt plus 0.3fil minus 3pt

{\bf AMS Classification numbers}\quad Primary:\quad \theprimaryclass

Secondary:\quad \thesecondaryclass

\vglue 5pt plus 0.3fil minus 2pt

{\bf Keywords:}\quad \thekeywords

\vglue 10pt plus 0.5fil minus 5pt

{\small  Proposed: \theproposer\hfill Received: \receiveddate\nl
Seconded: \theseconders\hfill 
\ifx\reviseddate\relax                         
Accepted: \accepteddate                        
\else
Revised: \reviseddate                          
\fi}
\eject
}       

\let\maketitlepage\maketitlep
\let\maketitle\maketitlepage

\font\phead=cmsl9 scaled 950
\font=cmsl9 scaled 1050
\font\pnum=cmbx10 scaled 913
\font\lnum=cmbx10 
\font\pfoot=cmsl9 scaled 950
\font=cmsl9 scaled 1050
\ifplaintex
\headline{\vbox to 0pt{\vskip -4.5mm\line{\small\phead\ifnum
\count0=\startpage ISSN 1364-0380 (on line)
1465-3060 (printed) \hfill {\pnum\folio}\else\ifodd\count0\def\\{ }%
\ifx\theshorttitle\relax\thetitle\else\theshorttitle\fi\hfill{\pnum\folio}
\else\def\\{ and }{\pnum\folio}\hfill\ifx\theshortauthors\relax\theauthors
\else\theshortauthors\fi\fi\fi}\vss}}
\footline{\vbox to 0pt{\vglue 0mm\line{\small\pfoot\ifnum\count0=\startpage
\copyright\ \gtp\hfill\else
\gt, Volume \thevolumenumber\ (\thevolumeyear)\hfill\fi}\vss
}}
\else
\makeatletter
\def\@oddhead{{\small\lhead\ifnum\count0=\startpage ISSN 1364-0380 (on line)
1465-3060 (printed) \hfill {\lnum\number\count0}\else\ifodd\count0
\def\\{ }\ifx\theshorttitle\relax \thetitle \else\theshorttitle\fi\hfill
{\lnum\number\count0}\else\def\\{ and }{\lnum\number\count0}
\hfill\ifx\theshortauthors\relax 
\theauthors\else\theshortauthors\fi\fi\fi}}\def\@evenhead{\@oddhead}
\def\@oddfoot{\small\lfoot\ifnum\count0=\startpage\copyright\ \gtp\hfill\else
\gt, Volume \thevolumenumber\ (\thevolumeyear)\hfill\fi}
\def\@evenfoot{\@oddfoot}
\makeatother
\fi

\newwrite\gtoutfile
\long\gdef\makeheadfile{  
{\def\\{, }\def\s{ }
\immediate\openout\gtoutfile head.xxx
\immediate\write\gtoutfile{Proxy-for: \ifx\theasciiauthors\relax
\theauthors\else\theasciiauthors\fi\s<\ifx\theasciiemail\relax\theemail\else\theasciiemail\fi>}
\immediate\write\gtoutfile{\noexpand\\}
\immediate\write\gtoutfile{Authors: \ifx\theasciiauthors\relax
\theauthors\else\theasciiauthors\fi}
{\def\\{ }\immediate\write\gtoutfile{Title: \ifx\theasciititle\relax
\thetitle\else\theasciititle\fi}}
\immediate\write\gtoutfile{Subj-class: GT or SG or MG etc}
\immediate\write\gtoutfile{MSC-class: \theprimaryclass\ifx\thesecondaryclass\relax\else, \thesecondaryclass\fi}
\immediate\write\gtoutfile{Journal-ref: Geom. Topol. \thevolumenumber
(\thevolumeyear) \startpage-\finishpage}
\immediate\write\gtoutfile{Comments: Published by Geometry and Topology at}
\immediate\write\gtoutfile{\s\s http://www.maths.warwick.ac.uk/gt/GTVol\thevolumenumber/paper\thepapernumber.abs.html}
\immediate\write\gtoutfile{\noexpand\\}
\immediate\write\gtoutfile{}
\ifx\theasciiabstract\relax
\immediate\write\gtoutfile{\theabstract}\else
\immediate\write\gtoutfile{\theasciiabstract}\fi
\immediate\write\gtoutfile{}
\immediate\write\gtoutfile{\noexpand\\}
\immediate\write\gtoutfile{}
\immediate\closeout\gtoutfile}}  

\def\maketitlepage{\maketitlep\makeheadfile}
\let\maketitle\maketitlepage

\lognumber{326}
\volumenumber{8}\papernumber{12}\volumeyear{2004}
\pagenumbers{511}{538}
\received{29 July 2003}
\revised{12 March 2004}
\published{13 March 2004}
\accepted{16 December 2004}
\proposed{Benson Farb}
\seconded{Jean-Pierre Otal, Shigeyuki Morita}

\usepackage{amsmath,amssymb}
\usepackage{graphicx,psfrag}

\newcommand\C{\mathbb C}

\newcommand\ZZ{{\mathbb Z}_2}
\newcommand\Z{\mathbb Z}
\newcommand\Q{\mathcal Q}
\newcommand\HA{\mathcal H}
\newcommand\R[1]{{\mathbb R}^{#1}}

\newcommand\D{d\hspace{-0.5pt}}
\newcommand\M{M^{^{^{\hspace{-6pt} \circ}}}}

\renewcommand{\mod}{\:\operatorname{mod}}

\DeclareMathOperator{\Ker}{Ker}
\DeclareMathOperator{\ind}{ind}

\newtheorem{Theorem}{Theorem}[section]
\newtheorem{Corollary}[Theorem]{Corollary}

\newtheorem*{Construction}{Construction}
\newtheorem*{NoNumberTheorem}{Theorem}
\newtheorem*{NoNumberProposition}{Proposition}

\theoremstyle{remark}
\newtheorem{Remark}{Remark}

\theoremstyle{definition}

\newtheorem*{Moduli:Spaces}{Moduli Spaces}
\newtheorem*{Stratification}{Stratification}
\newtheorem*{Teichmull}{Teichm\"uller Geodesic Flow}
\begin{document}
\title[Parity of the spin structure]
{Parity of the spin structure defined\\by a quadratic
differential}

\author{Erwan Lanneau}
\address{Institut de math\'ematiques de Luminy \\
Case 907, 163 Avenue de Luminy \\
F-13288 Marseille Cedex 9, France}
\email{lanneau@iml.univ-mrs.fr}


\asciikeywords{Quadratic differentials, Teichmueller
geodesic flow, moduli space, measured foliations, spin structure}

\begin{abstract}
According to the work of Kontsevich--Zorich,
the invariant that classifies non-hyperelliptic connected
components of the moduli spaces of Abelian differentials with
prescribed singularities, is the {\it parity of the spin
structure}.

We show that for the moduli space of quadratic differentials, the
spin structure is constant on every stratum where it is defined.
In particular this disproves the conjecture that it classifies the
non-hyperelliptic connected components of the strata of quadratic
differentials with prescribed singularities. An explicit formula
for the parity of the spin structure is given.
\end{abstract}

\asciiabstract{According to the work of Kontsevich-Zorich,
the invariant that classifies non-hyperelliptic connected
components of the moduli spaces of Abelian differentials with
prescribed singularities,is the parity of the spin
structure.

We show that for the moduli space of quadratic differentials, the
spin structure is constant on every stratum where it is defined.
In particular this disproves the conjecture that it classifies the
non-hyperelliptic connected components of the strata of quadratic
differentials with prescribed singularities. An explicit formula
for the parity of the spin structure is given.}

\maketitlepage

\section{Introduction}

Quadratic differentials and moduli spaces of quadratic
differentials are natural objects in Teichm\"uller theory. In
particular, they are related to so called interval exchange
transformations and billiards in rational polygons. The moduli
spaces of Abelian differentials (denoted $\HA_g$) and the moduli
spaces of quadratic differentials (denoted $\Q_g$) are naturally stratified by the genus of
surfaces and the type of the singularities of the differentials. Due to
works of Veech \cite{Veech:82} and Rauzy \cite{Rauzy}, the topology
of the strata of these spaces is related to the dynamic of these intervals
exchange transformations.

A fundamental result independently proved by Masur and by Veech
(1982) asserts that the Teichm\"uller geodesic flow acts
ergodically on each connected component of any stratum of the
moduli spaces $\HA_g$ and $\Q_g$. Recently, Kontsevich and Zorich
have described the set of components for a particular type of
strata: those of moduli space of Abelian differential $\HA_g$
(see~\cite{Kontsevich:Zorich}). In order
to obtain their classification, they use two invariants:
the hyperellipticity and the spin structure. In a previous
paper \cite{Lanneau:02:hyperelliptic}, we give three families
of non-connected strata of $\Q_g$, using the
hyperelliptic invariant. Moreover, we show that there are no other
hyperelliptic components of the strata of $\Q_g$.
See also~\cite{Lanneau:02:classification}
for a complete description of the set of components of the strata of the
moduli space of quadratic differentials $\Q_g$.

In this paper we are interested in the calculation of the second
invariant on the elements of $\Q_g$. Here we prove that it is constant on every stratum of
$\Q_g$ where it is defined. Thus the parity of the spin structure cannot distinguish 
two differentials with the same singularity pattern. It is therefore a negative answer to 
a question of Kontsevich and Zorich on the classification of the non-hyperelliptic 
connected component by this spin structure. In addition we give an explicit
formula for its calculation knowing the type of the singularities.
In the Appendix we recall the relationship between rational billiards
and quadratic differentials and then, using our formula, we give
some applications in terms of billiards.

\subsection{Background}

On a Riemann surface $M^2_g$ of genus $g$, a meromorphic quadratic differential
is locally defined by the
form $\psi=f(z)(\D z)^2$ where $z$ is a local coordinate. In this paper
we consider quadratic differentials having only simple poles, if any:
that is the functions $f$ are meromorphic with only simple poles, if any.

Here we will use the geometric point of view of quadratic
differentials: a flat structure with cone type singularities. Such
surfaces are those which possesses locally the geometry of a
standard cone. We can define it by a flat Riemannian metric with
specific isolated singularities. The standard cone possesses a
unique invariant: it is the angle at the vertex. Here we consider
only {\it half-translation} flat surfaces: parallel transport of a
tangent vector along any closed path either brings the vector
$\vec{v}$ back to itself or brings it to the centrally-symmetric
vector $-\vec{v}$. This implies that the cone angle at any
singularity of the metric is an integer multiple of $\pi$.

One can see that these Euclidian structures are induced by
quadratic differentials by the following way. Let $\psi$ be a
meromorphic quadratic differential on a Riemann surface $M^2_g$.
Then it is possible to choose a canonical atlas on $M \backslash
\{$singularities$\}$ such that $\psi=dz^2$ in any coordinate
chart. As $dz^2=dw^2$ implies $z=\pm w +const$, we see that the
charts of the atlas are identified either by a translation or by a
translation composed with a central symmetry. Thus a meromorphic
quadratic differential $\psi$ induces a half-translation flat
structure on $M \backslash \{$singularities$\}$. On a small chart
which contains a singularity, coordinate $z$ can be chosen in such
way that $\psi=z^kdz^2$, where $k$ is the order of the singularity
($k = 0$ corresponds to a regular point, $k=-1$ corresponds to a
pole and $k > 0$ corresponds to a true zero). It is easy to check
that in a neighborhood of a singularity of $\psi$, the metric has
a cone type singularity with the cone angle $(k+2) \pi$.

Conversely, a half-translation structure on a Riemann surface
$M$ and a choice of a distinguished ``vertical'' direction defines
a complex structure and a meromorphic quadratic differential
$\psi$ on $M$. In many cases it is very convenient to present a
quadratic differential with some specific properties by an
appropriate flat surface. Consider, for example, a polygon in the
complex plane $\C$ with the following property of the boundary:
the sides of the polygon are distributed into pairs, where the
sides in each pair are parallel and have equal length. Identifying
the corresponding sides of the boundary by translations and
central symmetries we obtain a Riemann surface with a natural
half-translation flat structure. The quadratic differential $dz^2$
on $\mathbb{C}$ gives a quadratic differential on this surface
with punctures. The punctures correspond to vertices of the
polygon; they produce the cone type singularities on the surface.
It is easy to see that the complex structure, and the quadratic
differential extends to these points, and that a singular point of
the flat metric with a cone angle $(k+2) \pi$ produces a
singularity of order $k$ (a pole if $k=-1$) of the quadratic
differential. (See also Figure~\ref{fig_1} which illustrates this
construction.)

\begin{Moduli:Spaces}
According to these definitions, one can define the moduli space
$\HA_g$ of {\it Abelian differentials} as the moduli space of
pairs $(M,\omega)$ where $\omega$ is a holomorphic $1-$form
defined on a Riemann surface $M$ of genus $g$. Here the term
moduli spaces means that the points $(M_1,\omega_1), \
(M_2,\omega_2)$ are identified if and only if there exists an
isomorphism $f \co  M_1 \rightarrow M_2$ (with respect to the complex structure)
affine in the canonical chart determined by $\omega_i$. In an equivalent way
we can ask that $f^\ast \omega_2=\omega_1$.

We can also define the moduli space $\Q_g$ of quadratic
differentials as the moduli space of pairs $(M,\psi)$ where $\psi$
is a meromorphic quadratic differential {\it which is not} the
global square of a $1-$form defined on $M$.
\end{Moduli:Spaces}

Recall that we denote by $(k_1,\dots,k_n)$ the orders of singularities of
$\psi$. A consequence of the Gauss-Bonet formula is that $\sum k_i = 4g-4$.

\begin{Stratification}
The moduli space $\Q_g$ is naturally stratified by the types of
singularities of the forms. We denote by $\Q(k_1,\dots,k_n)\subset \Q_g$ the
stratum of quadratic differentials $[M^2_g,\psi]\in\Q_g$ which are
not the squares of Abelian differentials, and which have the
singularity pattern $(k_1,\dots,k_n)$, where $k_i \in
\{-1,0,1,2,\dots\}$. We use the exponential
notation $k_i^m$ for $k_i,k_i,k_i,\dots,k_i$ repeated $m$ times. For
example, $\Q(1^4,8,2,3^2)$ stands for $\Q(1,1,1,1,8,2,3,3)$.

We can also consider the moduli space of Abelian differentials or quadratic
differentials which are the global square of Abelian differentials. We denote these
spaces by $\mathcal{H}_g$. And, if $\overrightarrow{k}$ is a vector
in $\mathbb{N}^n$ with $\sum k_i = 2g-2$, we denote by
$\mathcal{H}(k_1,\dots,k_n)$ the corresponding stratum.
\end{Stratification}

\begin{Teichmull}
The group SL$(2,\mathbb{R})$ acts on these spaces $\Q_g$ and
$\HA_g$ in the following way. For a matrix $A
\in$~SL$(2,\mathbb{R})$ and a point $[M,\psi]$ we define
$A[M;\psi]$ by the point $[M;A\ \psi]$. Notation $A\ \psi$ means
$A$ acts linearly in the canonical charts determined by $\psi$. By
definition, this action preserves each stratum. The following
one-parameter subgroup are of special interest:
$$
g_t = \left(
\begin{array}{cc}
e^{t/2} & 0 \\
0    & e^{-t/2} \
\end{array}  \right)
$$
We call it the Teichm\"uller geodesic flow.
\end{Teichmull}

In classical works, Masur and Veech discovered that the geodesic
flow on the moduli space of quadratic differentials $\HA_g \sqcup \Q_g$
is related to the theory of the so-called interval exchange transformation.
The phase space of the Teichm\"uller geodesic flow can be seen as the
cotangent bundle to the Teichm\"uller space, and it can be interpreted as the moduli space of
pairs consisting of a Riemann surface endowed with a holomorphic quadratic
differential. It is well known that the flow preserves the natural stratification and
that each stratum carries a complex algebraic orbifold structure. Moreover,
Masur and Smillie have proved that all of these strata, except four particular
cases in low genera, are non-empty.

The study of the topology of the strata comes from a fundamental Theorem,
independently proved by Masur and Veech.

\begin{NoNumberTheorem}[Masur, Veech, 1982]
The Teichm\"uller geodesic flow acts ergodically on every connected component
of any stratum with respect to a finite equivalent Lebesgue measure.
\end{NoNumberTheorem}

Kontsevich and Zorich have recently described the set of connected components for a
particular type of strata: one of the moduli space of {\it Abelian} differentials
$\HA_g$ \cite{Kontsevich:Zorich}.
In two papers \cite{Lanneau:02:hyperelliptic, Lanneau:02:classification}, we have obtained the ``complementary''
case, namely the description of components of the moduli space of {\it quadratic}
differentials $\Q_g$.

In their classification, Kontsevich and Zorich use two invariants
to obtain a complete description of the components: the hyperellipticity and the
parity of the spin structure. These invariants allow them to show that each stratum
of the moduli space $\HA_g$ possesses at most $3$ components.

In this paper, we are interested in the computation of the second invariant in $\Q_g$,
that is the parity of the spin-structure of an arbitrary quadratic differential.
Our main result is that this invariant is constant on every stratum of $\Q_g$ where
it is defined. It is therefore a negative result in contrast to that of Kontsevich--Zorich.

\subsection{Formulation of the statement}

Given an Abelian differential $\omega$ with even zeroes, one can associate
the following equation in the Picard group Pic$(M^2_g)$
\begin{equation*}
2\cdot K = D(\omega) = \sum_{i=1}^n 2\cdot k_i P_i
\end{equation*}
We define the spin structure determined by $\omega$ by
$\sum_{i=1}^n k_i P_i$. The parity of the spin structure
is defined as $\Phi(\omega) = \dim |\sum_{i=1}^n k_i P_i|$ mod $2$.
Following works of M~Atiyah and D~Mumford \cite{Atiyah, Mumford},
it can be shown that this integer is invariant under continuous deformation.
Moreover the spin-structure is the basic invariant which allows to classify all
non-hyperelliptic connected components of any stratum of $\HA_g$.

In section~\ref{s:covering:construction}, we consider a canonical local mapping
(induced by the standard orientating double covering)
$$
\Q(k_1,\dots,k_n) \longrightarrow \mathcal{H}(\tilde{k}_1,\dots,\tilde{k}_r).
$$
Using this mapping, we associate to each pair $[M^2_g,\psi]$ an
integer $0$ or $1$ which we defined as the parity of the spin structure
determined by $\psi$. This integer is invariant under continuous deformation
of the point $[M^2_g,\psi]$ in the given stratum.

In section~\ref{sec:invariance} we show
that for quadratic differentials the spin-structure is constant on
every stratum. Moreover we give an explicit formula to determine
it knowing the singularity pattern $(k_1,\dots,k_n)$ of the
stratum. This gives a negative response to a question of Kontsevich and Zorich
that it may distinguish the non-hyperelliptic connected components of some
strata $\Q(k_1,\dots,k_n)$. More precisely, we will show

\begin{Theorem}
Let $\psi$ be a quadratic differential on a Riemann surface such that
$\psi$ does not possess zeroes of order $k$ with $k=2$ mod $4$.
This condition assures that the spin structure of $\psi$ is well defined.
Then the parity of the spin structure of $\psi$ is
independent of the choice of $\psi$ in a stratum $\Q(k_1,\dots,k_l)$.
\end{Theorem}

In addition, we end section~\ref{sec:invariance} by the following description Theorem
in terms of $k_i$.

\begin{Theorem}
Let $\psi$ be a meromorphic quadratic differential on a Riemann
surface $M$ with singularity pattern $\Q(k_1,\dots,k_l)$. We also
require that $k_i\not = 2$ mod $4$.
Let $n_{+1}$ be the number of zeros of $\psi$ of
degrees $k_i = 1 \mod 4$, let $n_{-1}$ be the number
of zeros of $\psi$ of degrees $k_j = 3 \mod  4$, and
suppose that the degrees of all the remaining zeros satisfy $k_r
= 0 \mod  4$. Then the parity of the spin structure
defined by $\psi$ is given by
$$
\Phi(\psi)=\left[\cfrac{|n_{+1}-n_{-1}|}{4}\right]\
\mod 2
$$
where square brackets denote the integer part.
\end{Theorem}

\subsection{Deformation of flat surfaces}

\begin{figure}[ht!]
\begin{center}
\includegraphics[width=2.5cm]{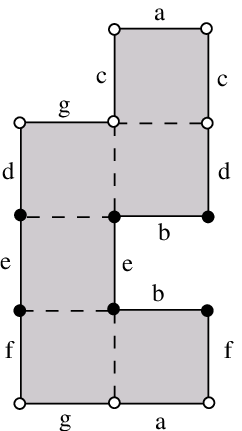}
\end{center}
\caption{
\label{fig_1}
A point in the stratum $\HA(2,2)$}
\end{figure}

Here, we explicitly present an example to illustrate that the problem of
classification of component is quite difficult. For this, let us consider
the domain given by Figure~\ref{fig_1} in $\mathbb{C}$ endowed with its complex structure
induced by the form $\D z$. We identify the opposite sides by a translation (with respect
to the numbered of sides). We get a Riemann surface of genus $g=3$ with an Abelian
differential induced from the form $\D z$. It is easy to see that, according to
previous notations, we obtain a point in the stratum $\HA(2,2)$.

\begin{figure}[ht!]\small
\begin{center}
\psfrag{m1}{$(M,\omega_1)$}  \psfrag{m2}{$(M,\omega_2)$}
\includegraphics[width=6.5cm]{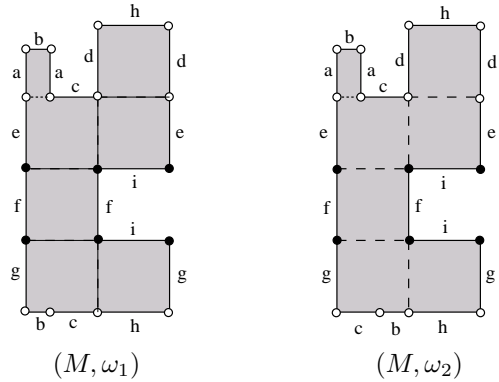}
\end{center}
\caption{
\label{fig_2}
Two points $(M,\omega_1)$ and $(M,\omega_2)$ in the stratum $\HA(2,4)$}
\end{figure}

Now we can make some surgery on this polygon: we glue a small rectangle on
the boundary of the polygon and identified the vertical boundary of the small
additional rectangle with appropriate translations (see Figure~\ref{fig_2} for details).
As above,  we make identifications of opposite sides according the numbered of sides.
We get a Riemann surface. It is not difficult to see that it has genus
$g=4$. The two Abelian differentials induced from the form $\D z$
in the complex plane have one zero of order $2$ and one zero of
order $4$ (given respectively by black and white bullets). Thus, with our above notation,
we obtain {\it two} points in the stratum $\HA(2,4)$. Let denote them respectively by
$(M,\omega_1)$ and $(M,\omega_2)$.

It is difficult to see directly if these two points belong in the
same connected component of the stratum $\HA(2,4)$ or not. In this
case, using the parity of the spin structure, we can check that
$\omega_1$ and $\omega_2$ are not in the same connected component
(one can see this by the fact that there exists an additional handle
and then we can apply Lemma~$11$ of~\cite{Kontsevich:Zorich}).

\begin{figure}[ht!]\small
\begin{center}
\psfrag{m3}{$(M,\omega_3)$}  \psfrag{m4}{$(M,\omega_4)$}
\includegraphics[width=6.5cm]{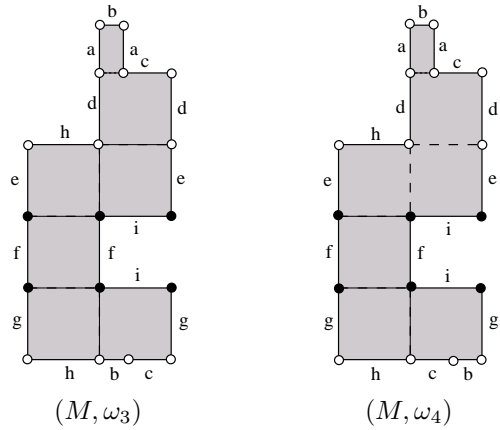}
\end{center}
\caption{
\label{fig_3}
Two points $(M,\omega_3)$ and $(M,\omega_4)$ in the stratum $\HA(2,4)$}
\end{figure}

However, we can perform an analogous surgery to the first one (see Figure~\ref{fig_3})

In the same way, we obtain two other points in the stratum $\HA(2,4)$. Let us denote
them by $[M,\omega_3]$ and $[M,\omega_4]$. Using the parity of the spin structure,
we can deduce that these two points {\it are} in the same component. In fact,
using notations of~\cite{Kontsevich:Zorich}, we have
$$
[M,\omega_1],\ [M,\omega_3],\ [M,\omega_4] \in \HA^{odd}(2,4)
\qquad \textrm{and} \qquad [M,\omega_2] \in \HA^{even}(2,4)
$$
However, it seems very difficult to prove this fact directly.

The paper has the following structure: we present some general
facts about quadratic differentials in
section~\ref{s:covering:construction}. In
section~\ref{sec:invariance}, we announce and prove the main
Theorem. In section~\ref{subsec:formula} we give the formula to
compute the parity of the spin structure. In the Appendix we recall
some classical constructions related to rational polygonal
billiards and we use our formula to obtain applications in this context.

\subsection{Acknowledgements}

I would like to thank Anton Zorich for discussion on this subject
and remarks concerning this paper. I thank the Institut de
Math\'ematiques de Luminy and the Max-Planck-Institut f\"ur
Mathematik at Bonn for excellent welcome and working conditions.

\section{Mapping of the moduli spaces induced by a ramified
covering of a fixed combinatorial type}
\label{s:covering:construction}

In this section we present some general information concerning the
moduli spaces of quadratic differentials. The proofs and the
details can be found in papers \cite{Masur:82, Veech:82, Masur:Smillie, Veech,
Kontsevich, Kontsevich:Zorich, Douady:Hubbard}.

\begin{NoNumberTheorem}[Masur and Smillie]
Consider a vector $(k_1,\dots,k_n)$ with all $k_i \in \mathbb{N}
\cup \{-1 \}$. Suppose that $\sum k_i = 0 \mod 4$ and $\sum k_i
\geq -4$. The corresponding stratum $\Q(k_1,\dots,k_n)$ is
non-empty with the following four exceptions:
$$
 \Q(\emptyset), \Q(1,-1)\ (\text{in genus } g=1)
 \quad \text{and}\quad
 \Q(4), \Q(1,3)\ (\text{in genus } g=2)
$$
\end{NoNumberTheorem}

\begin{NoNumberTheorem}[Masur, Veech]
Any stratum $Q(k_1,\dots,k_n)$ is a complex orbifold of dimension

\centerline {${\rm dim}_{\mathbb{C}} \Q(k_1,\dots,k_n) = 2g + n -2$.}
\end{NoNumberTheorem}

\begin{NoNumberProposition}[Kontsevich]
Any stratum $\Q(k_1,\dots,k_n)$ with $\sum k_i = -4$ is connected.
\end{NoNumberProposition}

Here we recall a classical construction of the ``orientating map''
of a quadratic differential. We refer to~\cite{Douady:Hubbard} for
a general reference.

\begin{Construction}[Canonical double covering]
\label{constr:canonical:2:covering}
Let $M^2_g$ be a Riemann surface and let $\psi$ be a quadratic
differential on it which is not a square of an Abelian
differential. There exists a canonical (ramified) double covering
$\pi \co  \tilde{M}^2_{\tilde{g}} \rightarrow M^2_g$ such that
$\pi^{*} \psi=\tilde{\omega}^2$, where $\tilde{\omega}$ is an
Abelian differential on $\tilde{M}^2_{\tilde{g}}$.

The images $P\in M^2_g$ of ramification points of the covering
$\pi$ are exactly the singularities of odd degrees of $\psi$. The
covering $\pi\co  \tilde{M}^2_{\tilde{g}} \rightarrow M^2_g$ is the
minimal (ramified) covering such that the quadratic differential
$\pi^\ast\psi$ becomes the square of an Abelian differential on
$\tilde{M}^2_{\tilde{g}}$.
\end{Construction}

\begin{proof} Consider an atlas
$(U_i,z_i)_i$ on ${\M}^2_g=M^2_g \backslash
\{\text{singularities of }\psi\}$ where we punctured all zeros
and poles of $\psi$. We assume that all the charts $U_i$ are
connected and simply-connected. The quadratic differential $\psi$
can be represented in this atlas by a collection of holomorphic
functions $f_i(z_i)$, where $z_i\in U_i$, satisfying the
relations:
$$
 f_i\left(z_i(z_j)\right)\cdot\left(\cfrac{\D z_i}{\D z_j}\right)^2 =
 f_j(z_j)\ \text{ on }U_i \cap U_j
$$
Since we have punctured all singularities of $\psi$ any function
$f_i(z_i)$ is nonzero at $U_i$. Consider two copies $U^{\pm}_i$ of
every chart $U_i$: one copy for every of the two branches
$g_i^\pm(z_i)$ of $g^\pm(z_i):=\sqrt{f_i(z_i)}$ (of course, the
assignment of ``$+$'' or ``$-$'' is not canonical). Now for every
$i$ identify the part of $U^{+}_i$ corresponding to $U_i\cap U_j$
with the part of one of $U^{\pm}_j$ corresponding to $U_j\cap U_i$
such that on the overlap the branches match:
$$
 g^{+}_i(z_i(z_j))\cdot\cfrac{\D z_i}{\D z_j} =
 g^{\pm}_j(z_j)\ \text{ on }U^{+}_i \cap U^{\pm}_j
$$
Apply the analogous identification to every $U^{-}_i$. We get a
Riemann surface with punctures provided with a holomorphic 1-form
$\tilde{\omega}$ on it, where $\tilde{\omega}$ is presented by the
collection of holomorphic functions $g^{\pm}_i$ in the local
charts. It is an easy exercise to check that filling the punctures
we get a closed Riemann surface $\tilde{M}^2_{\tilde{g}}$, and
that $\tilde{\omega}$ extends to an Abelian differential on it. We
get a canonical (possibly ramified) double
covering $\pi\co \tilde{M}^2_{\tilde{g}}\to M^2_g$ such that
$\pi^{\ast}\psi=\tilde{\omega}^2$.

By construction the only points of the base $M^2_g$ where the
covering might be ramified are the singularities of $\psi$. In a
small neighborhood of zero of even degree $2k$ of $\psi$ we can
choose coordinates in which $\psi$ is presented as $z^{2k} (\D
z)^2$. In this chart we get two distinct branches $\pm z^k \D z$
of the square root. Thus the zeros of even degrees of $\psi$ and
the marked points are the regular points of the covering $\pi$.
However, it easy to see that the covering $\pi$ has exactly a
ramification point over any zero of odd degree and over any simple
pole of $\psi$.
\end{proof}

\begin{Remark}
We can consider this construction for any quadratic differential. Note
that corresponding quadratic differential is the square of an Abelian
differential if and only if $\tilde{M_g}$ is non-connected.
\end{Remark}

Now, consider a non-degenerate point $[M,\varphi] \in
\Q(k_1,\dots,k_n)$, that is a non-orbifoldic point. Deforming
slightly the initial point $[M^2_g,\psi_0]\in\Q(k_1,\dots,k_n)$ we
can consider the canonical double ramified covering over the
deformed Riemann surface of the same combinatorial type as the
covering $\pi$. This new covering has exactly the same relation
between the positions and types of the ramification points and the
degrees and position of singularities of the deformed quadratic
differential. This means that the induced quadratic differential
$\pi^{\ast}\psi$ has the same singularity pattern. Choosing one of
two branches of $\omega=\pm \sqrt{\pi^\ast \psi}$, we get a local
mapping:
\begin{gather*}
\Q(k_1, \dots, k_n) \to \mathcal{H}(\tilde{k}_1,\dots,\tilde{k}_m)\\
[M^2_g,\psi]\mapsto (\tilde{M}^2_g,\pi^{\ast}\psi)
\end{gather*}
\begin{Remark}
A detailed proof of this fact can be found in \cite{Lanneau:02:hyperelliptic}.
\end{Remark}

\section{Parity of a spin structure defined by a quadratic differential}
\label{s:parity:of:the:spin:structure}

\subsection{Spin structure defined by an Abelian differential}
We first recall the algebraic-geometric definition of the {\it
spin structure} given by an Abelian differential, see
M~Atiyah~\cite{Atiyah}; see also~\cite{Kontsevich:Zorich}; see
D~Johnson~\cite{Johnson} for a topological definition.

\subsubsection{Definition of a spin structure}
Let $\omega$ be a $1$-form with only even singularities.
There are $2^{2g}$ solutions of the equation $2 D = K(\omega)$ in the
divisor group where $K(\omega)$ is the canonical divisor determined
by $\omega$. A spin structure is the choice of $D$ in the Picard group of
the surface.
For an Abelian differential with only even zeros, one can write
$$
K(\omega) = 2k_1P_1 + \dots + 2k_nP_n
$$
With these notation, we declare that the spin structure defined by the form
$\omega$ on the complex curve $M$ is just the divisor
$D=k_1P_1 + \dots + k_nP_n$. Thus a point $(M,\omega)$ gives
{\it canonically} a spin structure.

\subsubsection{Parity of a spin structure}
The dimension of the linear space $|D|$ may have quite
different values for different choice of $D$. For example,
in genus $1$, the dimension of the spaces given by the three non-null
solutions have non-zero dimension. We declare that the dimension modulo
$2$ of this linear space is the parity of the spin structure $D$ and we denote it
by $\Phi(D)$. On a curve of genus $g \geq 1$, M~Atiyah~\cite{Atiyah}
proved that there are $2^{g-1}(2^g+1)$ odd spin structure and
$2^{2g}-2^{g-1}(2^g+1)$ even spin structure.

\subsubsection{Calculation of the parity of a spin structure}
Let $M$ be a Riemann surface with an Abelian differential
$\omega$. Denote by $\M$ the surface $M$ punctured at the
singularities of $\omega$. Consider the flat metric on $\M$
defined by $\omega$. The corresponding holonomy representation in
the linear group is trivial: a parallel transport of a tangent
vector along any closed loop brings the vector to itself. Thus we
get a canonical trivialization of the tangent bundle
$T_\ast\M$ to the punctured Riemann surface $\M$: the
tangent space at any point of $\M$ is canonically identified
with a sample Euclidean plane $\R{2}$.

Let $\gamma$ be a closed smooth oriented connected curve on $M$
avoiding singularities of $\omega$. Using the trivialization of
the tangent bundle to $\M$ we can construct the Gauss map
($\mathcal{G}\co  \gamma \to S^1$) as follows: we associate to a
point $x\in\gamma$ the image of the normalized tangent vector to
$\gamma$ at $x$ under the trivialization map $T_1 \M \to
S^1\subset\R{2}$. Here the unit tangent bundle $T_1\M\subset
T_\ast\M$ to $\M$ is defined in terms of the flat metric on $\M$
defined by $\omega$.

The {\it index} of the closed curve $\gamma$ in the flat metric
defined by $\omega$ is defined to be the degree of the Gauss
map
$$
\ind_\gamma(\omega):=\deg(\mathcal{G}\co\gamma\to S^1).
$$
In other words, following the closed curve $\gamma$ in the
positive direction and measuring how the tangent vector turns in
the flat structure we get the angle $\ind_\gamma(\omega)\cdot
2\pi$ as a total turn along $\gamma$. For example, given a zero of
order $k \geq 0$ and a small curve $\gamma$ around it, we get for
the index of $\gamma$ $(k+1)2\pi$.

\begin{Remark}
\label{rk:index}
If $\triangle$ is a small disc in $M$ avoiding singularities, by
definition $[\partial \triangle]=0$ in $H_1 (\M; \mathbb Z_2)$.
Thus we have
\begin{multline*}
\Omega_\omega(\partial \triangle)=\Omega_\omega(\partial \triangle
+
\partial \triangle) =
 \\
 = \Omega_\omega(\partial \triangle) + \Omega_\omega(\partial \triangle)
+ \partial \triangle \circ \partial \triangle =
\\
= 2 \cdot \Omega_\omega(\partial \triangle) + 0 = 0  \textrm{ mod
} 2
\end{multline*}
where $\circ$ stands for the standard intersection form on
$H_1(M;\mathbb Z_2)$ mod $2$.

Moreover, if we have only one singularity inside $\triangle$ of
order $k$ we have
$$
\Omega_\omega(\partial \triangle) = \textrm{ ind}_{\partial
\triangle}(\omega) + 1  =  k \textrm{ mod } 2.
$$
We obtain $k = 0 \textrm{ mod } 2$ thus $\omega$ must possesses
only singularities of even order.
\end{Remark}

Remark~\ref{rk:index} implies that when the Abelian differential
$\omega$ has zeros of even degree only, the residue
$\ind_\gamma(\omega)\mod 2$ depends only on the homology class
$[\gamma]\in H_1(M;\ZZ)$ of the smooth closed connected path$\gamma$.

This invariance allows to define the function $\Omega_\omega\co 
H_1(M;\ZZ)\to \ZZ$: to evaluate $\Omega_\omega$ on a cycle $c\in
H_1(M;\ZZ)$ we represent $c\in H_1(M;\ZZ)$ by a closed smooth
oriented connected curve $\gamma$ on $M$ avoiding singularities of
$\omega$, and define
$$
\Omega_\omega(c):=\ind_\gamma(\omega)+1\mod 2.
$$
The function $ \Omega_\omega $ is, actually, a {\it quadratic
form} on $H_1(M;\ZZ)$, that is
\begin{equation}
\label{eq:Arf}
   \Omega_\omega(c_1+c_2)=
   \Omega_\omega(c_1)+\Omega_\omega(c_2)+c_1\circ c_2.
\end{equation}
Let $\{a_1,b_1,\dots,a_g,b_g\}$ be a symplectic basis of
$H_1(M,\Z)$. We define the {\it Arf-invariant} of the
quadratic form $\Omega_\omega$ as follows:
\begin{equation*}
\Phi(\Omega_\omega)=
\sum_{i=1}^{g} \Omega_\omega(a_i)\Omega_\omega(b_i)
\mod 2
\end{equation*}
According to \cite{Arf}, the number defined by the above formula,
is independent of the choice of a symplectic basis of cycles.
In addition, Johnson \cite{Johnson} proved that the set of
quadratic forms and spin structures is one to one. Finally, Johnson
has shown that the two numbers (modulo $2$) defined by the parity of the spin structure
and the Arf invariant of the corresponding quadratic form must coincide. Thus we have
$$
\Phi(\omega)=\Phi(K(\omega))=\textrm{dim} | D(\omega) | \ \textrm{ mod } 2
= \Phi(\Omega_\omega)
$$
where $D(\omega)=k_1P_1+\cdots+k_nP_n$, and $2k_i$ is the multiplicity of
the zero $P_i$ of $\omega$.

\subsubsection{Spin structure on a deformed curve}
According the result of M~Atiyah~\cite{Atiyah} and D~Mumford~\cite{Mumford},
the dimension of the linear space $ | D(\omega) | $ modulo $2$ is invariant under
continuous deformations of the Abelian differential $\omega$
inside the corresponding stratum, and hence, it is constant
for every connected component of any stratum, where it is defined.

\subsection{Spin structure on a double covering
defined by a quadratic differential}

Now we want to define the parity of the spin structure for a
quadratic differential. Recall that this number is well defined
only for Abelian differentials with zeroes of even orders.

Let $M$ be a Riemann surface and $\psi$ a quadratic differential
on $M$ which is not the square of any Abelian differential. We
suppose that orders of all singularities of $\psi$ are different
from $2$ modulo $4$. We consider the canonical (branched) double
covering $\pi \co  \tilde{M} \rightarrow M$ over $M$ such that
$\pi^\ast \psi=\omega^2$ (see Construction in
Section~\ref{s:covering:construction}). Then we can check that
condition on the degree of the zeroes of $\psi$ implies that
$\omega$ possesses only zeros of even order (Lemma~$1$
in~\cite{Lanneau:02:hyperelliptic}). So we can apply the notion of
parity of the spin structure to $\omega$.
With these notation, we declare that the parity of spin structure of
$\psi$ is the parity of the spin structure determined by $\omega$:
\begin{equation*}
\Phi(\psi) \stackrel{\textrm{Def}}{:=} \Phi(\omega)
\end{equation*}
As $\omega$ and $-\omega$ are in the same connected component (consider
the path $e^{i \pi t}\omega$ for $t \in [0,1]$),
we have $\Phi(\omega)=\Phi(-\omega)$. So $\Phi(\psi)$ is well-defined.

\section{Invariance of the parity of the spin structure}
\label{sec:invariance}

In this Section, we prove the announced result.

\subsection{Monodromy representation for the canonical double covering}

Let $M$ be a Riemann surface with $\psi$ a quadratic
differential, that is not the square of an Abelian
differential. Let $\pi\co \tilde{M}\to M$ be the canonical
(ramified) double covering such that the pullback $\pi^\ast\psi$
of $\psi$ to the covering surface $\tilde{M}$ becomes the global square of
an Abelian differential,
$\pi^\ast(\psi)=\omega^2$. Removing from $\tilde{M}$ the ramification points
(if there are any), and removing the singularities of $\psi$ of odd
degrees (if there are any) from $M$ we obtain a regular
$\ZZ$-covering
$$
\pi\co  \tilde{M}\setminus\{\text{ramification points}\}\to M\setminus
\{\text{singularities of }\psi \text{ of odd degree}\}
$$
Denote $\M:=M\setminus \{\text{singularities of }\psi \text{ of
odd degree}\}$. We get a monodromy representation $m\co  \pi_1(\M)
\to \ZZ$. Since $\ZZ$ is an Abelian group the monodromy
representation quotients as $\pi_1(\M) \to H_1(\M;\Z{})\to \ZZ$,
and even as $\pi_1(\M) \to H_1(\M;\ZZ)\to \ZZ$. We denote the
monodromy homomorphism $m\colon  H_1(\M;\ZZ)\to \ZZ$ by the same symbol
$m$.

A closed connected path $\gamma$ on $\M$ lifts to a closed
path on the covering surface $\tilde{M}$ if and only if the
homology class $[\gamma]\in H_1(\M;\ZZ)$ belongs to
$\Ker(m)$. Note, however, that if $\psi$ does have singularities
of odd degree it is always possible to modify a closed connected
path $\gamma$ on $\M$ in such way that the modified closed
connected path $\gamma'$ stays in the same homology class
$[\gamma]=[\gamma']\in H_1(M;\ZZ)$, but the monodromy along
$\gamma$ and $\gamma'$ are different, $m(\gamma)\neq
m(\gamma')$. The reason is that the
homology classes of $\gamma$ and $\gamma'$ in homology of the
{\it punctured} surface $H_1(\M;\ZZ)$
are different. In example presented by Figure~\ref{fig:cycles}.
the difference corresponds to a cycle turning around the
puncture, and such a cycle has nontrivial monodromy.

Note, that if a cycle in $H_1(\M;\ZZ)$ is represented by a
smooth closed connected curve $\gamma$ on $M$ which does not pass
through singularities of $\psi$, the monodromy $m(\gamma)$ can be
calculated geometrically. It is trivial if and
only if the holonomy of the flat metric defined by $\psi$ along
this path is trivial --- a parallel transport of a vector along the
path brings the vector to itself.

\subsection{Main result}

Now we have all necessary tools to prove the following Theorem.

\begin{Theorem}
\label{theo:spin}
The parity of the spin structure of a quadratic differential
$\psi$, that is not the square of an Abelian differential, is
independent of the choice of $\psi$ in a stratum
$\Q(k_1,\dots,k_l)$.
\end{Theorem}

\begin{proof}[Proof of the Theorem~\ref{theo:spin}]
We first present an idea of the proof. Let $M$ be a Riemann
surface with a quadratic differential $\psi$ on it. Let $\pi \co 
\tilde{M} \rightarrow M$ be the canonical double covering such
that $\pi^\ast\psi=\omega^2$ ; by hypothesis, $\tilde{M}$ is
connected. Let $2n$ be the number of singularities of $\psi$ of
odd degree; the double covering $\pi \co  \tilde{M} \rightarrow M$
has the ramification points exactly over these points. By the
Riemann--Hurwitz formula the genus of the covering surface
$\tilde{M}$ equals $2g+n-1$. Hence
\begin{equation}
\label{eq:dimH}
 \dim H_1(\tilde{M};\mathbb{Z}_2) = 4g+2n-2.
\end{equation}
We start the proof with a construction of a special basis of
cycles on $\tilde{M}$:
$$
a_1^+, b_1^+, a_1^-, b_1^-, \dots,
a_{g}^+, b_{g}^+, a_{g}^-, b_{g}^-,
\tilde{c}_1, \dots, \tilde{c}_{n-2}
$$
The basis is partially orthogonalized: the cycles $a_i^\pm,
b_i^\pm$ are obtained by lifting to $\tilde{M}$ closed paths on
$M$ representing a symplectic basis in $H_1(M;\ZZ)$. Thus
$a_i^+\circ b_i^+=a_i^-\circ b_i^-=1$, and the other intersection
indices for these cycles are zero. The cycles $\tilde{c}_j$
correspond to the ``cuts'' joining the consecutive pairs $P_j,
P_{j+1}$ of ramification points of the covering; these cycles are
completely analogous to the cycles from a standard Riemann--Hurwitz
basis on a hyperelliptic surface.

The basis is constructed in such way that
\begin{equation}
\label{eq:abpm}
\ind_{a^+_i}\omega=\ind_{a^{-}_i}\omega \qquad
\ind_{b_i^+}\omega=\ind_{b_i^-}\omega,
\end{equation}
and $\ind_{\tilde{c}_j}\omega$ is easily expressed in terms of
degrees $k_j$ and $k_{j+1}$ of the corresponding zeros $P_j,
P_{j+1}$ of $\psi$. Thus:
\begin{multline}
\label{eq:impact:c}
 \Phi(\omega)=
\sum_{i=1}^{2g}(\ind_{a_i^+}(\omega)+1)(\ind_{b_i^+}(\omega)+1) +
\sum_{i=1}^{2g}(\ind_{a_i^-}(\omega)+1)(\ind_{b_i^-}(\omega)+1) +
 \\
+\text{ impact of the cycles } \tilde{c}_j\ \mod  2 =
\\
= \text{ impact of the cycles } \tilde{c}_j\ \mod  2
\end{multline}
The remaining part of the theorem is an exercise in linear algebra
and in arithmetic. We present now the complete proof.

We first treat the special case when $\psi$ has no singularities
of odd degree at all. Note that in this case the covering
$\pi\co \tilde{M}\to M$ is a regular double covering. In particular,
the monodromy of the covering along a closed path
depends only on the homology class of the
path. A closed path $\gamma$ on $M$
lifts to a closed path $\tilde{\gamma}$ on $\tilde{M}$ if and only
if the corresponding cycle $[\gamma]$ belongs to $\Ker(m)$,
where $m\co  H_1(M;\ZZ)\to \ZZ$ is the monodromy homomorphism.

In the case when $\psi$ has no singularities of odd degree the
basis of cycles on $\tilde{M}$ which we are going to construct
has the form
$$
a_1^+, b_1^+, a_1^-, b_1^-, \dots, a_{g-1}^+, b_{g-1}^+,
a_{g-1}^-, b_{g-1}^-,\ a_g^+, \tilde{b}_g \ ,
$$
where there are only two cycles $a_g^+, \tilde{b}_g$
corresponding to index $g$. We obtain this basis of cycles on the
covering surface $\tilde{M}$ using a special basis of
cycles on the underlying surface $M$.

Let $a_1,b_1,\dots,a_g,b_g$ be a canonical basis of
$H_1(M,\mathbb{Z}_2)$, where $a_i\circ b_i=1$ for $i=1, \dots, g$,
and all the other intersection indices are trivial. Since the
quadratic differential $\psi$ is not a global square of an Abelian
differential, there exists a cycle $c\in H_1(M,\mathbb{Z}_2)$ with
nontrivial monodromy. Without loss of generality we may assume that
this cycle is $b_g$, that is $m(b_g)=1$.

Let us perform now a change of the basis in $H_1(M,\mathbb{Z}_2)$
to obtain a canonical basis in $H_1(M;\mathbb{Z}_2)$ with the
additional property:
$$
\ind_{a_i}\psi=\ind_{b_i}\psi=0  \text{ for } i=1,\dots,g-1
\quad\text{ and }\quad
\ind_{a_g}\psi=0, \ \ind_{b_g}\psi=1
$$
We replace those $a_i$, $i<g$, and those $b_j$, $j<g$, that have
nontrivial monodromy, by $a_i+b_g$ and $b_j+b_g$ correspondingly.
Now we have $m(a_i)=m(b_i)=0  \text{ for } i=1,\dots,g-1$. It
remains to adjust $a_g$ to get a canonical basis with desired
properties.

We now use this basis on the underlying surface $M$ to construct a
basis of cycles on the covering surface $\tilde{M}$. Since
$m(a_i)=m(b_i)=0  \text{ for } i=1,\dots,g-1$, a
closed path on $M$ representing any cycle from this collection can
be lifted to a closed path on the covering surface $\tilde{M}$ in
two different ways; the two corresponding lifts are disjoint.
In this way we obtain the representatives of the cycles $a_1^+,
b_1^+, a_1^-, b_1^-, \dots, a_{g-1}^+, b_{g-1}^+$ on $\tilde{M}$.
We may assign the superscript indices $\pm$ in such a way that
$a^+_i\circ b^+_i=a^-_i\circ b^-_i=1$ and all the other
intersection indices are zero. We now take a closed path on $M$
representing the cycle $a_g$ and we choose one of the two lifts of
this path. The resulting closed path on $\tilde{M}$ gives us a
cycle $a^+_g\in H_1(\tilde{M},\mathbb{Z}_2)$. Since
$m(b_g)=1$ we need to take a double lift of a
representative of this cycle to get a  closed path on $\tilde{M}$;
it gives us $\tilde{b}_g$. We get a collection of $4g-2$ cycles in
$H_1(\tilde{M},\mathbb{Z}_2)$. By construction the matrix of the
intersection form for this collection has the canonical form; in
particular, it is non-degenerate. Since the number of cycles equals
the dimension of the homology group $H_1(\tilde{M};\mathbb{Z}_2)$,
see equation~\eqref{eq:dimH}, we get a canonical basis of cycles
on the covering surface $\tilde{M}$.

Let us compute the parity of the spin structure of the Abelian
differential $\omega$ on $\tilde{M}$ using this basis.
By construction the basis is canonical, and
equations~\eqref{eq:abpm} are valid for any $i=1,\dots,g-1$.
Thus:
\begin{multline*}
\Phi(\omega)=
\sum_{i=1}^{g-1}(\ind_{a_i^+}\omega+1)(\ind_{b_i^+}\omega+1) +
\sum_{i=1}^{g-1}(\ind_{a_i^-}\omega+1)(\ind_{b_i^-}\omega+1) +
\\
+ (\ind_{a_g^+}(\omega)+1)(\ind_{\tilde{b}_g}(\omega)+1)
 \mod  2 =
\\
=(\ind_{a_g^+}(\omega)+1)(\ind_{\tilde{b}_g}(\omega)+1)
 \mod  2
\end{multline*}
By assumption $m(b_g)=1$. This means that the holonomy of the flat
structure corresponding to $\psi$ along a smooth path on $M$
representing $b_g$ is an odd multiple of the angle $\pi$. Hence
the holonomy of the flat structure corresponding to $\omega$ along
a smooth path on $\tilde{M}$ which is a double cover of the
corresponding path on $M$ is an odd multiple of the angle $2\pi$.
This means that $\ind_{\tilde{b}_g}\omega = 1 \mod  2$, which
implies $\Phi(\omega)=0$. We proved the Theorem in a special case
when $\psi$ does not have poles no zeroes of odd degree. Moreover,
we showed that the parity of the spin structure is {\it always}
even in this special case.

\begin{figure}[ht!]\small
\begin{center}
\psfrag{pj}{$P_j$}   \psfrag{pi+j}{$P_{j+1}$}
\psfrag{g}{$\gamma$}   \psfrag{g'}{$\gamma'$}
\psfrag{cj}{$c_j$}   \psfrag{p}{$P_{2n}$}
\psfrag{m}{$M^2_g$}
\includegraphics[width=8cm]{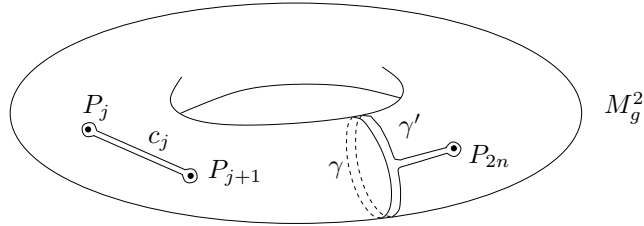}
\end{center}
\caption{
\label{fig:cycles}
The points $P_j, P_{j+1}$, and $P_{2n}$ of the underlying surface $M$
correspond to singularities of $\psi$ of odd degree. The
paths $\gamma$ and $\gamma'$ are homologous in $H_1(M;\ZZ)$,
but holonomy along $\gamma$ and along $\gamma'$ is different.
Holonomy along the path $c_j$ is obviously trivial.
}
\end{figure}

Suppose now that $\psi$ has singularities of odd degree. Note that
the number of such singularities is always even, we denote it by
$2n$. Now the covering $\pi\co \tilde{M}\to M$ is ramified, and the
holonomy homeomorphism is well-defined only for the punctured
surface $M\setminus \{\text{ramification points}\}$. Adding to a
path $\gamma$ an appendix which goes to the ramification point,
turns around it and then comes back, we get a new path
representing the same homology class in $H_1(M; \ZZ)$ but having
different holonomy, see Figure~\ref{fig:cycles}. In particular,
now we can choose a collection of smooth  closed paths with trivial
holonomy on the underlying surface $M$ representing a canonical
basis of cycles $a_i, b_i$ in $H_1(M; \ZZ)$. We denote these
paths by the same symbols $a_i, b_i$ as the corresponding cycles.

Using  this collection of paths on the underlying surface $M$ we
construct a special basis of cycles
$$
a_1^+, b_1^+, a_1^-, b_1^-, \dots, a_{g}^+, b_{g}^+, a_{g}^-,
b_{g}^-, c_1, \dots, c_{n-2}
$$
on the covering surface $\tilde{M}$. Here the cycles $a^\pm_i,
b^\pm_i$ correspond to the lifts of the paths $a_i, b_i$. We
assign the superscript indices $\pm$ to these lifts in such way
that $a^+_i\circ b^+_i=a^-_i\circ b^-_i=1$.

The cycles $c_j$ are analogous to the ones from a standard
Riemann--Hurwitz basis on a hyperelliptic surface. They are
constructed as follows. For any $j=1,\dots,2n-2$ consider a path
on the underlying surface $M$ joining $P_j$ and $P_{j+1}$. We may
choose this path in such way that it does not intersect any
of the paths $a_i, b_i$, and that the corresponding ``broken line''
$P_1,P_2, \dots, P_{2g-1}$ does not have self-intersections. Deform
these paths as on the Figure~\ref{fig:cycles} to obtain closed
loops $c_j$ with trivial monodromy. Lifting $c_j$ to a closed path
$\tilde{c}_j$ on $\tilde{M}$ (in any of two possible ways) we get
the desired cycles in $H_1(\tilde{M};\ZZ)$.

The subspace spanned by the cycles
$a^+_i,a^-_i,b^+_i,b^-_i$, where $i=1,\dots,g$, is orthogonal to
the subspace spanned by the cycles $c_j$, where $j=1,\dots,2n-2$.
In the same way of the standard Riemann--Hurwitz basis on a hyperelliptic 
surface, we also have $\tilde{c}_j\circ\tilde{c}_{j+1}=1$ for
$j=1,\dots,2n-3$, and $\tilde{c}_{j_1}\circ\tilde{c}_{j_2}=0$ if
$|j_2-j_1|>1$. Thus the  matrix of the intersection form for this
collection of cycles is non-degenerate. According to
equation~\eqref{eq:dimH} the number $4g+2n-2$
of cycles equals the dimension of the homology group
$H_1(\tilde{M};\ZZ)$. Hence we
get a basis of cycles on the covering surface $\tilde{M}$.

Let us compute the parity of the spin structure of the Abelian
differential $\omega$ on $\tilde{M}$ using this basis. The
subspace $V_1$ spanned by the cycles $a^+_i,a^-_i,b^+_i,b^-_i$,
where $i=1,\dots,g$, is orthogonal to the subspace spanned by the
cycles $c_j$, where $j=1,\dots,2n-2$. By construction the cycles
$a^+_i,b^+_i,a^-_i,b^-_i$ give us a canonical basis of the
subspace $V_1$, and for any $i=1,\dots,g$ we have
relations~\eqref{eq:abpm}. Hence
$$
\sum_{i=1}^{2g}(\ind_{a_i^+}\omega+1)(\ind_{b_i^+}\omega+1) =
\sum_{i=1}^{2g}(\ind_{a_i^-}\omega+1)(\ind_{b_i^-}\omega+1)
$$
which justifies equation~\eqref{eq:impact:c}. \\
We have shown that the parity of the spin structure depends only of
computation of index of cycles $\tilde{c}_j$. Thus this number
depends only of degree of the singularities of $\omega$ that
is the number $k_i$. This completes the proof of Theorem~\ref{theo:spin}.
\end{proof}

\subsection{Explicit formula}
\label{subsec:formula}

In this section, we give an explicit formula for the parity of the spin
structure.

\begin{Theorem}
\label{theo:formula}
Let $\psi$ be a meromorphic quadratic differential on a Riemann
surface $M$ with singularity pattern $\Q(k_1,\dots,k_l)$.
Let $n_{+1}$ be the number of zeros of $\psi$ of
degrees $k_i = 1 \mod 4$, let $n_{-1}$ be the number
of zeros of $\psi$ of degrees $k_j = 3 \mod  4$, and
suppose that the degrees of all the remaining zeros satisfy $k_r
= 0 \mod  4$. Then the parity of the spin structure
defined by $\psi$ is given by
$$
\Phi(\psi)=\left[\cfrac{|n_{+1}-n_{-1}|}{4}\right]\
\mod 2
$$
where square brackets denote the integer part.
\end{Theorem}

\begin{proof}[Proof of the Theorem~\ref{theo:formula}]
Let $M$ be a Riemann surface endowed with a quadratic
differential $\psi$. Let $\pi \co  \tilde{M} \rightarrow M$ be
canonical double covering such that $\pi^\ast \psi = \omega^2$.
We want to compute $\Phi(\psi)=\Phi(\omega)$.
If $\psi$ has no singularities of odd degree at all (that is neither
zeroes of odd degree, nor poles) then we already showed
that $\Phi=0$ which proves the Theorem~\ref{theo:formula} in a
special case. \\
Suppose now that $\psi$ has singularities of odd degree.
We use the notations in the proof of the Theorem~\ref{theo:spin}.
According to the equation~\eqref{eq:impact:c}, we have:
\begin{multline*}
 \Phi(\psi)=
\sum_{i=1}^{2g}(\ind_{a_i^+}(\omega)+1)(\ind_{b_i^+}(\omega)+1) +
\sum_{i=1}^{2g}(\ind_{a_i^-}(\omega)+1)(\ind_{b_i^-}(\omega)+1) +
 \\
+\text{ impact of the cycles } \tilde{c}_j\ \mod  2 =
\\
= \text{ impact of the cycles } \tilde{c}_j\ \mod  2
\end{multline*}
If $\psi$ has only two singularities of odd degree, our basis does not
have any cycles $\tilde{c}_j$ at all, so in this case the Theorem is proved.
Suppose now that $\psi$ has more than two singularities of odd
degree.

Let us compute now $\ind_{\tilde{c}_j}\omega$. We have to compute
how the tangent vector to the corresponding path $\tilde{c}_j$ on
$\tilde{M}$ turns in the flat structure defined by $\omega$, see
Figure~\ref{fig:index:c}. The counterclockwise direction is chosen
as a direction of the positive turn.

\begin{figure}[ht!]\small
\begin{center}
\psfrag{p1}{$\tilde{P_1}$}  \psfrag{p2}{$\tilde{P_2}$}
\psfrag{a1}{$A_1$}  \psfrag{a5}{$A_5$}
\psfrag{a2}{$A_2$}  \psfrag{a6}{$A_6$}
\psfrag{a3}{$A_3$}  \psfrag{a7}{$A_7$}
\psfrag{a4}{$A_4$}  \psfrag{a8}{$A_8$}

\includegraphics[width=8cm]{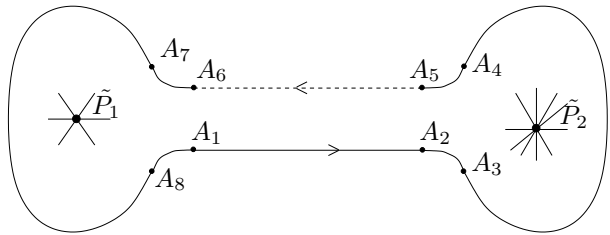}
\end{center}
\caption{
\label{fig:index:c}
A path $\tilde{c}_j$ on the covering surface obtained by lifting the path
$\tilde{c}_j$. Since $\tilde{P}_j, \tilde{P}_{j+1}\in \tilde{M}$ are ramification
points, the segments $A_1 A_2$ and $A_6 A_5$ are located on ``different sheets''
of the covering $\tilde{M}\to M$, and are not near to one another.
}
\end{figure}

We get some angle $\phi$ following the part $A_1 A_2$ of the path
which goes from one singularity to another. Then we make a turn by
$-\pi/2$ going along $A_2 A_3$. Turning around the singularity
(the path $A_3 A_4$) we get the angle $(k_{j+1}+2)\pi$, which is
followed by another turn by $-\pi/2$ now along $A_4 A_5$. The path
$A_5 A_6$ gives the turn by $-\phi$, which is followed by another
$-\pi/2$ along $A_6 A_7$. Turning around singularity along $A_7
A_8$ we get the angle $(k_j+2)\pi$, and the loop is completed by the
path $A_8 A_1$ giving one more turn by $-\pi/2$. All together this
gives $(k_j+k_{j+1}+2)\pi=((k_j+k_{j+1})/2+1)\cdot 2\pi$. So we obtain
$\ind_{\tilde{c}_j}(\omega)=(k_j+k_{j+1})/2+1$.
Thus the value of the $\ZZ$-quadratic form $\Omega_\omega$ on the
cycle $\tilde{c}_j$ is equal to
$$
\Omega_\omega(\tilde{c}_j)=\ind_{\tilde{c}_j}(\omega)+1 = \cfrac{k_j
+ k_{j+1}}{2} \mod 2
$$
Applying the Gram--Schmidt algorithm to the family $\tilde{c}_j$
we get the following symplectic basis $\alpha_j$, $\beta_j$,
$j=1,\dots,n-1$ in the subspace $V_2\subset H_1(\tilde{M};\ZZ)$
spanned by the cycles $\tilde{c}_j$:
$$
\begin{array}{cc}
\alpha_1=\tilde{c}_1                        & \beta_1=\tilde{c}_2          \\
\alpha_2=\tilde{c}_3+\alpha_1               & \beta_2=\tilde{c}_4          \\
\alpha_3=\tilde{c}_5+\alpha_2               & \beta_3=\tilde{c}_6          \\
\vdots                                      &\vdots                    \\
\alpha_{n-1}=\tilde{c}_{2n-3}+\alpha_{n-2}  & \beta_{n-1}=\tilde{c}_{2n-2}
\end{array}
$$
Then
$$
\Phi(\psi):=\Phi(\omega)=\sum_{j=1}^{n-1} \Omega_\omega(\alpha_j)
\cdot \Omega_\omega(\beta_j) \mod  2.
$$
We use formula~\eqref{eq:Arf} to evaluate the $\ZZ$-quadratic form
$\Omega_\omega$ for this new bases:
$$
\Omega_\omega(\alpha_j)=
\Omega_\omega(\tilde{c}_1+\tilde{c}_3+\dots+\tilde{c}_{2j-1})=
\Omega_\omega(\tilde{c}_1)+\Omega_\omega(\tilde{c}_3)+
\dots+\Omega_\omega(\tilde{c}_{2j-1}),
$$
where we take into consideration that $\tilde{c}_{2j_1-1} \circ \tilde{c}_{2j_2-1} = 0$
for any $j_1,j_2$. Since
$\Omega_\omega(\beta_i)= \Omega_\omega(\tilde{c}_{2i})$ we
finally obtain:
\begin{equation}
\label{eq:phi:long}
\Phi(\psi)=
\begin{cases}
\cfrac{1}{4} \sum_{j=1}^{n-1} (k_1+k_2+\dots+k_{2j})(k_{2j} +
k_{2j+1})
\mod 2 & \text{ if } n \geq 2\\
0 & \text{ if } n=0,1
\end{cases}
\end{equation}
To complete the proof of the Theorem we need to apply several
elementary arithmetic arguments. Assume that $n>1$.
First note that all the integers $k_j$ are odd. Hence inside any pair of parentheses
we have an even number which implies that we may
replace the numbers $k_j$ in~\eqref{eq:phi:long} by their
residues modulo $4$ without affecting the total sum modulo $2$.

We may enumerate the singularities of odd degrees in such way
that the first $2m$ ones $P_1, \dots, P_{2m}$ the residues $k_j \mod 4$
alternate, and for the remaining $2(n-m-1)$ ones
$P_{2m+1}, \dots, P_{2n-2}$ the residues
$k_i \mod 4$ are the same. Then, for any $j\le m$ we get
$k_1+\dots+k_{2j}= 0 \mod 4$, and hence:
\begin{multline*}
\Phi(\psi)=
\cfrac{1}{4} \sum_{j=1}^{n-1} (k_1+k_2+\dots+k_{2j})(k_{2j} +
k_{2j+1})
\mod 2 = \\
= \cfrac{1}{4} \smash{\sum_{j=m+1}^{n-1}}\vrule width 0pt depth 18pt (k_{2m+1}+\dots+k_{2j})(k_{2j} +
k_{2j+1}) \mod 2
\end{multline*}
Thus it is sufficient to check the formula in the statement
of the Theorem only for two cases: when all zeros have degree $1$,
and when all zeros have degree $3$.
Taking into consideration that for the total sum we have
$$
k_1+\dots+k_{2n-2}=0 \mod 4
$$
we easily obtain the desired relation.
\end{proof}

\subsection{Ergodic components of the Teichm\"uller geodesic flow}

Recall that the two invariants which classify connected components of the strata
of $\HA_g$ are the hyperellipticity and the spin structure.
In the previous section, we have shown that the parity of the spin structure is constant
on each stratum of the moduli space of meromorphic quadratic differential $\Q_g$.
Moreover, we can calculate all hyperelliptic components of the moduli space of quadratic
differentials. In~\cite{Lanneau:02:hyperelliptic}, we proved that
the strata $\Q(12)$ and $\Q(-1,9)$ do not possess an hyperelliptic component.
However, it was proved by A~Zorich by a direct computation of the
corresponding {\it extended Rauzy classes} that each of these two strata has exactly
two distinct connected components. Using formula of Theorem~\ref{theo:formula},
we get that the spin structure is even on these two strata and so, it does not distinguish
the component; which gives a negative response to a question of Kontsevich and Zorich.

\appendix

\section{Some applications}
\label{application:billiard}

\subsection{Billiard Flow}
Let $P$ denote a polygon in $\mathbb{R}^2$. The billiard flow is given by the motion
of a point with the usual optical reflection rule on the boundary $\partial P$ of $P$.
A line element of this flow (the geodesic flow) is given by a point $x \in P$ and a
direction $\theta \in \mathbb{R}/2\pi\mathbb{Z}$ on the unit tangent bundle of $P$. Orbits
of this flow fail to have continuations when they hit the boundary of $P$.  We would
like trajectories to reflect off the boundary. If $e_i$ is an edge of $P$ and
$\rho_i \co  S^1 \rightarrow S^1$ represent the reflexion through $e_i$ then we identify
$(p,v)$ with $(p,\rho_i(v))$ for each $p \in e_i$.

Let $\Gamma \subset O(2)$ be the group generated by the reflections in sides. We are
interested only by the case when $\Gamma$ is finite; in such case the polygon is
called rational. An equivalent condition, when $P$ is simply connected, is that all
angles are rational multiple of $\pi$. For a rational billiard, a classical construction
(see~\cite{Masur:Tabachnikov} for a nice review of rational billiards)
gives rise to a flat surface $\tilde{M}^2_g$ of genus $g$ i.e. a Riemann surface
endowed with an Abelian differential. We recall here the construction.

Let $P \subset \C$ and $\Gamma$ as above. Take $|\Gamma|$ disjoint
copies of $P$, each rotated by an element of $\Gamma$. For each
copy $P_c$ of $P$ and each reflection $r \in \Gamma$, glue each
edge $E_c$ of $P_c$ to the edge $r(E_c)$ of $r(P_c)$. When the
group $\Gamma$  is finite, the result is a compact Riemann surface
$\tilde{M}$
$$
\tilde{M}= \left. \left( \bigsqcup_{\gamma \in \Gamma} \gamma(P)
\right) \ \right/ \  \sim
$$
where $\sim$ is the relation above. The form $\D z$ on each copy $\gamma(P)$, $\gamma \in
\Gamma$, induces a holomorphic $1-$form $\tilde{\omega}$ on $\tilde{M}$.
It is easy to check that the singularities of $\tilde \omega$ are
located at the vertices of the copies of $P$. Moreover, the billiard flow
descends to the geodesic flow on the flat surface $(\tilde{M},\tilde{\omega})$. \\
When the billiard $P$ has symmetries, we can construct the smaller surface
denoted by $(\tilde{M_1},\tilde{\omega_1})$, obtained by identifying $\gamma_1(P)$ and
$\gamma_2(P)$ if they differ by a translation. We can deduce
$(\tilde{M},\tilde{\omega})$ from $(\tilde{M_1},\tilde{\omega_1})$ by
a finite covering
$$
(\tilde{M},\tilde{\omega}) \rightarrow (\tilde{M_1},\tilde{\omega_1}).
$$
The genus of $\tilde{M}$ is given in terms of the angles of $P$.

Thus a rational billiard defines a point $[\tilde{M},\tilde
\omega]$ in some stratum of moduli space $\mathcal{H}_g$. It is
easy to compute the singularity pattern of $(\tilde{M},\tilde
\omega)$ knowing the polygon $P$ and hence identify the stratum.
It is, somehow, more complicated to identify the connected
component of the stratum. The invariance of the spin-structure on
each stratum of quadratic differential gives information about the
connected components of the stratum which contains the surface
$\tilde{M}$.

\subsection{Billiard table and Abelian differential}
For this section we refer to~\cite{Masur:Tabachnikov}. Let $P$ be
a {\it rational} polygonal billiard: that is a rational polygon
with the billiard flow inside. Let $(m_i / n_i) \cdot \pi$ be the
angles of $\pi$ with the convention that gcd$(n_i,m_i)=1$. Let
$N=$lcm$(n_i)$ be the least common multiple of $n_i$. Denote by
$k$ the number of sides of $P$. Let $\tilde{M}$ be the Riemann
surface with the holomorphic $1-$form $\tilde{\omega}$ which
arises from $P$. We can calculate the genus $\tilde{g}$ of
$\tilde{M}$ as follows:
$$
\tilde{g}=1+\cfrac{N}{2}(k-2-\sum \cfrac{1}{n_i})
$$
We can also calculate singularities of $\tilde{\omega}$. Each
vertex of $P$ of angle $(m_i / n_i) \cdot \pi$, induces $N/n_i$ singularities of
order $m_i-1$ of the form $\tilde{\omega}$. If $m_i=1$, the singularity is a
{\it fake} zero (regular point of $\tilde \omega$) and we do not mark it.

\subsection{Billiard table and quadratic differential}

\subsubsection{Construction}

Here we present an alternative construction which associates a
surface $M$ with a {\it quadratic} differential $\psi$ from a
rational billiard $P$ (see~\cite{Fox:Kershner}). 
The image of the holonomy map for the flat
metric defined by $\psi$ on $M$ is contained in $\{ \pm
\textrm{Id} \}$. Let $\mathbb{P}^1$ be the double of $P$: we
consider two copies of $P$; one a top of the other and we identify
the boundaries of these two copies gluing the corresponding edges.
We obtain the topological sphere with a Euclidean metric. The
segments coming from edges of the polygon are non-singular for
this metric, so it has only isolated singularities coming from the
vertices of the polygon. The singularities are conical
singularities of angle rational multiple of $\pi$. We construct
the standard covering of this surface to obtain a translation
surface (not necessarily with non-trivial holonomy). The quadratic
differential $\psi$ which we obtain is a square of an Abelian
differential if and only if the $N=lcm(n_i)$ is odd. Otherwise the
surfaces $(\tilde{M},\tilde \omega)$ and $(M,\psi)$ are related by
the standard double covering.
$$
(\tilde{M},\tilde \omega) \rightarrow (M,\psi)
$$
When $N$ is odd, the two constructions coincide. We compute the
singularity pattern of the surface $(M,\psi)$ as follow. Let
$p_i/q_i$ denote the rational factor of $\pi/2$ for angles of $P$
with $p_i$ and $q_i$ coprime. Let $Q$ be the {\it lcm} of $q_i$.
Then for each vertex, the quadratic differential $\psi$ on $M$
possesses $Q/q_i$ singularities of order $p_i-2$.

\subsubsection{Formula}

When the invariant surface $(\tilde M,\tilde \omega)$ constructed
by a rational billiard belongs to a non-connected stratum, we can
apply Theorem~\ref{theo:formula} to identify the connected
component which contain the point $[\tilde{M},\tilde{\omega}]$.
Let $P$ be a rational billiard. Let the angles of $P$ have the
form $(m_i / n_i) \pi$ with $m_i=2k_i+1$. This condition
guarantees that the spin structure of the corresponding Abelian
differential is well defined; that is all zeroes of $\omega$ are
even. We enumerate the angles as follows: for $i=1,\dots,r_1$,
let $(m_i / n_i) \pi$ be the angles with $n_i$ even and $k_i$
even; for $i=r_1+1,\dots,r_1+r_2$, let $(m_i / n_i) \pi$ be the
angles with $n_i$ even and $k_i$ odd; finally for
$i=1+r_1+r_2,\dots,r+r_1+r_2$, let $(m_i / n_i) \pi$ be the angles
with $n_i$ odd. We denote by $N$ the lcm of the $n_i$. Then we
have the following
\begin{Corollary}
\label{cor:formula} Let $P$ be a billiard table as above.
Then the Abelian differential $\tilde \omega$ corresponding
to $P$ has the following parity of the spin-structure:
$$
\Phi(\tilde{\omega})=\left[\cfrac{N}{4} \ \cdot \ \left|
\sum_{i=1}^{r_1} \cfrac{1}{n_i} - \sum_{i=1+r_1}^{r_1+r_2}
\cfrac{1}{n_i}  \right| \ \right]\ \mod 2
$$
\end{Corollary}

\proof
First we calculate angles with the form $\cfrac{p_i}{q_i} \cdot
\cfrac{\pi}{2}$ with $p_i$ and $q_i$ coprimes. With considerations
above, it is not difficult to see that we can take $(p_i;q_i) =
(m_i;\cfrac{n_i}{2})$ for $i=1,\dots,r_1+r_2$ and $(p_i;q_i) = (2
\cdot m_i;n_i)$ for $i=1+r_1+r_2,\dots,r+r_1+r_2$. The degrees of
the singularities of the corresponding quadratic differential are
given by the $p_i$. For $i=1,\dots,r_1+r_2$, we obtain zeroes of
degree $p_i-2=m_i-2=2k_i-1$; for $i=1+r_1+r_2,\dots,r+r_1+r_2$, we
obtain zeroes of degree $p_i-2=2 \cdot (m_i-1)$. In the last case,
all zeroes are even so it does not appear in the formula of
Theorem~\ref{theo:formula}. When $k_i$ is even, the corresponding
$p_i$ is equal to $-1$ modulo $4$ and $k_i$ is odd, the
corresponding $p_i$ is equal to $+1$ modulo $4$. Let us calculate
the multiplicities of the zeroes. Let $Q$ denote the lcm of the
$q_i$. Now, the formula of the corollary follows from
Theorem~\ref{theo:formula} and the following fact
$$
Q=\cfrac{N}{2} \ ; \qquad  q_i = \cfrac{n_i}{2}\eqno{\qed}
$$

Quite often the translation surface $(\tilde M,\tilde \omega)$
corresponding to a rational billiard $P$ is hyperelliptic.
Corollary~\ref{cor:formula} allows one to determine the case when
it is not, see~\ref{appendix:hyper} . Though the surfaces $(\tilde
M,\tilde \omega)$ arising from billiards form subsets of measure
zero in corresponding strata, the properties of the ambient strata
have strong influence on the dynamics of the billiard. In
particular billiards inducing hyperelliptic and non-hyperelliptic
surfaces have different behavior.

\subsubsection{Hyperelliptic Components}
\label{appendix:hyper}

For moduli space of Abelian differentials, the parity of the spin
structure of the hyperelliptic components is well defined.
In~\cite{Kontsevich:Zorich}, Kontsevich and Zorich give a formula
for calculation of these numbers for {\it any} hyperelliptic
component. This is given by the following Proposition

\begin{NoNumberProposition}[Kontsevich,~Zorich]
The parity of the spin structure for a translation surface in a
hyperelliptic component with one singularity and for a translation surface in
a hyperelliptic component with two singularities is given, respectively, by
\begin{equation}
\label{formul:1}
\Phi(\HA(2g-2)^{hyp})=\left[ \cfrac{g+1}{2} \right] \mod 2
\end{equation}
and
\begin{equation}
\label{formul:2}
\Phi(\HA(g-1,g-1)^{hyp})=\cfrac{g+1}{2} \mod 2, \textrm{ for $g$ odd.}
\end{equation}
\end{NoNumberProposition}

\medskip

\noindent {\bf Example}\qua
Consider the polygon given by the triangle
$(11\pi/14;\pi/7;\pi/14)$. We consider the billiard flow inside
this polygon. The number $N$ appearing in the classical
construction $A.2.$ is $N=14$. In particular it is even. The
corresponding Abelian differential $(\tilde{M},\tilde \omega)$ has
$2$ fake zeroes and a zero of order $10$. Thus it determines a
point in the stratum $\mathcal{H}(10)$. The surface $\tilde{M}$
has genus $6$. We can calculate the singularity pattern of the
corresponding quadratic differential. The angles are given by the
rational factor of $\pi$: $\left(\cfrac{2\cdot5
+1}{14} \ ; \ \cfrac{2\cdot0+1}{14} \ ; \ \cfrac{1}{7} \right)$. Thus by
Corollary~\ref{cor:formula}, the parity of the spin structure of
$(\tilde{M},\tilde \omega)$ is even.

According to formula~\ref{formul:1}, the parity
of the spin structure of the hyperelliptic component of the stratum
$\mathcal{H}(10)$ is odd and so the surface $\tilde{M}$ is not
hyperelliptic, and $[\tilde M,\tilde \omega] \in
\mathcal{H}^{\textrm{odd}}(10)$.

\end{document}